\documentclass[a4paper, 12pt]{article}
\usepackage{amsmath,amsthm, amssymb, esint}
\usepackage{graphicx}
\usepackage{graphics}
\usepackage{marginnote}
\usepackage{framed}
\usepackage{color}
\definecolor{red1}{rgb}{1.00,0.00,0.00}

\newcommand{\Addresses}{{
  \bigskip
  \footnotesize

  \textsc{Department of Mathematics, Saarland University, P.O. Box 151150,  Saar- br{\"u}cken 66041, Germany} and
  \textsc{ Faculty of Mathematics and Mechanics, St. Petersburg State University, Universitetskii pr. 28,  St. Petersburg 198504, Russia}\par\nopagebreak
  \textit{E-mail address:} \texttt{darya@math.uni-sb.de}

  \bigskip

  \textsc{Faculty of Mathematics and Mechanics, St. Petersburg State University, Universitetskii pr. 28,  St. Petersburg 198504, Russia}\par\nopagebreak
  \textit{E-mail address:} \texttt{uraltsev@pdmi.ras.ru}
}}

\newtheorem{theorem}{\bf Theorem}[section]

\newtheorem{remark}{\bf Remark}[section]
\newtheorem{definition}{\bf Definition}[section]
\title{Free boundaries in problems with hysteresis
}
\author{D.E.\,Apushkinskaya and N.N.\,Uraltseva} 

\begin{document}
\maketitle

\begin{abstract}
In this note we present a survey concerning parabolic free boundary problems involving a discontinuous hysteresis operator.  Such problems describe biological and chemical processes "with memory" in which various substances interact according to hysteresis law. 

Our main objective is to discuss the structure of the free boundaries and the properties of the so-called "strong solutions" belonging to the anisotropic Sobolev class $W^{2,1}_q$ with sufficiently large~$q$. Several open problems in this direction are proposed as well.
\end{abstract}

\section{Introduction}

The paper concerns with parabolic equations containing discontinuos hysteresis operator on the right-hand side. For simplicity we restrict our consideration to the most basic discontinuous hysteresis operator - the so-called non-ideal relay. 

\begin{figure}[!h]
\centering\includegraphics[width=4in]{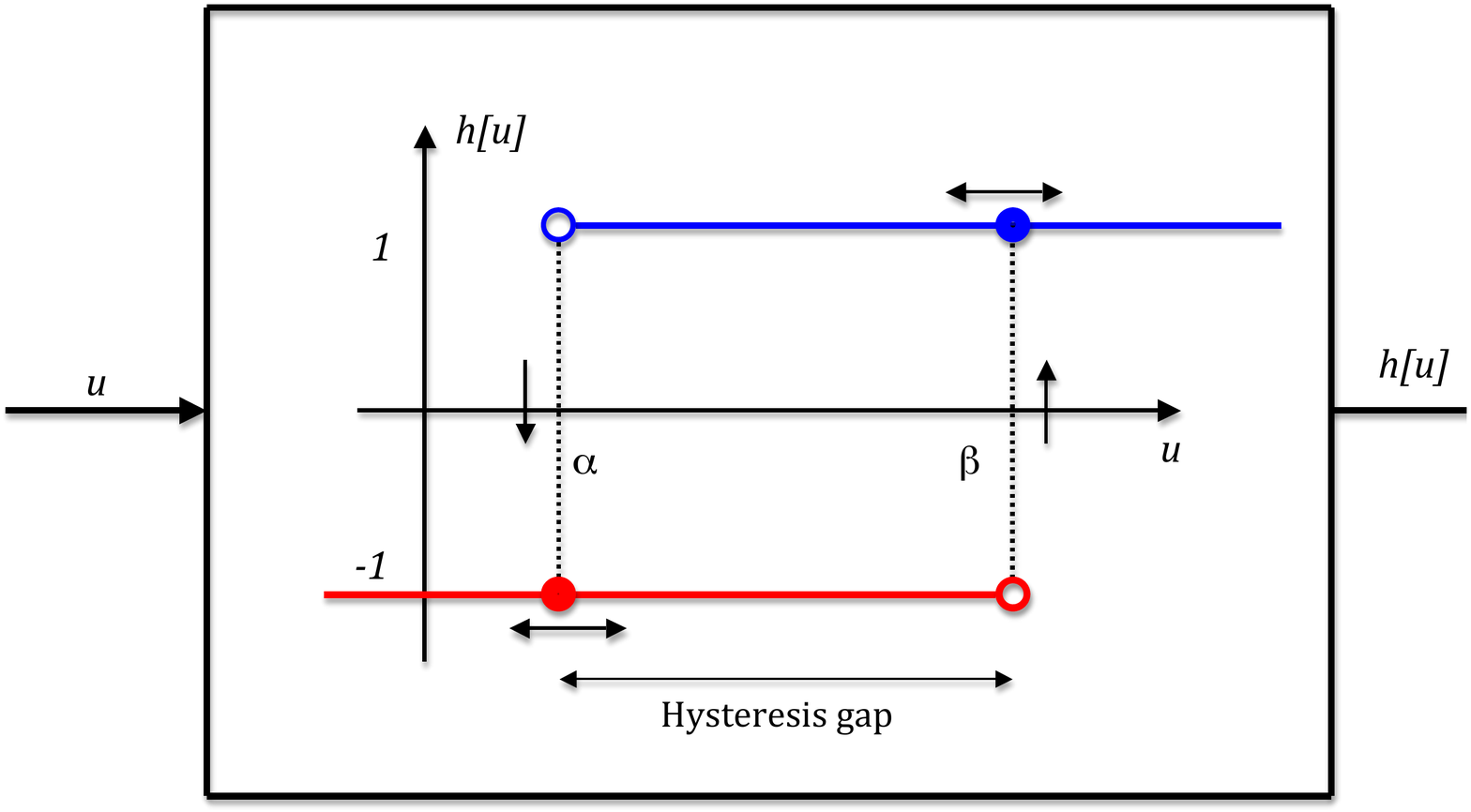}
\caption{Non-ideal relay}
\label{fig_1}
\end{figure}

If the input function $u$ is less than a lower treshold $\alpha \in \mathbb{R}$, then the output $h[u]$ of the non-ideal relay is equal $-1$. While $u$ increasing, the output remains equal $-1$ until the input  reaches an upper treshold $\beta \in \mathbb{R}$ - at this moment the output switches from $-1$ to $1$. Further increasing of $u$ does not change the output state. Observe that $h[u]$ switches back to the state $-1$ when the input $u$ afterwards decreases to  $\alpha$. This behaviour  is illustrated in Figure~\ref{fig_1}. It is evident that the non-ideal relay operator takes the path of a rectangular loop and its next state depends on its past state.

Examples of parabolic equations with non-ideal relay arise in various biological, technological and chemical processes 
(see, for instance, \cite{Vis86, Alt85, HJ80,HJP84,K06}, and references therein).

\subsection{Statement of the problem}
We study solutions of the nonlinear  parabolic equation
\begin{equation} \label{main-equation}
\Delta u -\partial_tu=h[u] \quad \text{in}\quad Q=\mathcal{U}\times (0,T],
\end{equation}
satisfying the initial condition
\begin{equation} \label{initial-data}
u\big|_{t=0}=\varphi (x), \qquad x\in \mathcal{U}.
\end{equation}
We also assume that $u$ satisfies either the Dirichlet or the Neumann boundary condition on the lateral surface of cylinder $Q$, i.e.,
\begin{equation} \label{boundary-conditions}
\text{either}\quad u=\psi_1 \quad \text{or} \quad \frac{\partial u}{\partial \textbf{n}}=\psi_2 \qquad \text{on} \quad \partial'Q:=\partial\mathcal{U} \times (0,T).
\end{equation}
Here $\Delta$ is the Laplace operator, $\mathcal{U}$ is a domain in $\mathbb{R}^n$ and $\partial\mathcal{U}$ is its boundary, $\varphi$ and $\psi_1$ (or $\psi_2$) are given functions,  while $h[u]$ stands for a non-ideal relay operator acting from $C(\mathcal{U}~\times~[0,T])$ to $\left\lbrace \pm 1\right\rbrace $. 

In order to define the operator $h[u]$ we  fix two numbers $\alpha$ and $\beta$ ($\alpha <\beta$) and consider a \textbf{multivalued} function
$$
f(s)=\left\lbrace 
\begin{aligned}
-&1,  \quad \text{for}\quad s\in ]-\infty, \alpha],\\
&1,  \quad \text{for}\quad s\in [\beta, +\infty[,\\
-&1\ \text{or}\ 1,  \quad \text{for} \quad s \in ]\alpha, \beta[.
\end{aligned}
\right.
$$
Assuming 
$
\varphi \in C^{1,1}(\mathcal{U})
$,
we suppose that the values of $h_0(x):=f(\varphi (x))$  are prescribed. We set 
\begin{equation}\label{initial-hysteresis}
h[u](x,0)=h_0(x), \qquad  x\in \mathcal{U}.
\end{equation}

After that for every point $z=(x,t)\in Q$ and for $u\in C(\overline{Q})$ the corresponding value of $h[u](z)$ is uniquely defined in the following manner. Let us denote by $E$ a set of points
$$
E:=\left\lbrace z\in Q : u(z)\leqslant \alpha\right\rbrace \cup \left\lbrace z\in Q : u(z) \geqslant \beta\right\rbrace
\cup \left\lbrace \mathcal{U} \times \left\lbrace 0\right\rbrace \right\rbrace  .
$$
In other words, $E$ is a set where $f(u(z))$ is well-defined.

If $z\in E$ then $h[u](z)=f(u(z))$. Otherwise, for $z=(x,t)\in Q$ such that $\alpha <u(z)<\beta$ we set
$$
h[u](x,t)=h[u](x, \hat{t} (x)),
$$
where
$$
\hat{t} (x)=\max\left\lbrace s : (x,s)\in E;\ s \leqslant t\right\rbrace .
$$

Observe that for fixed $x$ a jump of $h[u](x,\cdot )$ can happen only on thresholds $\left\lbrace u(x, \cdot ) =\alpha\right\rbrace $ and $\left\lbrace u(x,\cdot )=\beta\right\rbrace $. Moreover, \textbf{"jump down"} (from $h=1$ to $h=-1$)  is possble on  $\left\lbrace u(x,\cdot )=\alpha\right\rbrace $ only, whereas \textbf{"jump up"} (from $h=-1$ to $h=1$)  is possible on $\left\lbrace u(x,\cdot )=\beta\right\rbrace$ only. 

\begin{remark}
It should be emphasized that the above definition of $h[u]$ excludes the case $\alpha=\beta$ from the consideration.
\end{remark}

\begin{definition} \label{def-strong-solutions}
We say that $u$ is a (strong) solution of  (\ref{main-equation})-(\ref{initial-hysteresis}) if 
\begin{itemize}
\item[(i)] $u\in W^{2,1}_q(Q)$, $q>n+2$, and $h[u]$ is generated by the function $u$ in 
view~of~(\ref{initial-hysteresis}) as described above,
\item[(ii)] $u$ satisfies (\ref{main-equation}) a.e. in $Q$,
\item[(iii)] $u$ satisfies (\ref{initial-data})-(\ref{boundary-conditions}) is the sense of traces.
\end{itemize}
\end{definition}

\begin{remark}
Recall that $W^{2,1}_{q}\left( Q\right) $ is the anisotropic Sobolev space with the norm
$$
\|u\|_{W^{2,1}_{q}(Q)}=\|\partial_tu\|_{q,Q}+\|D^2u\|_{q,Q}+\|u\|_{q,Q},
$$
where $\|\cdot \|_{q,Q}$ denotes the norm in $L^q(Q)$ with $1<q<\infty$.
\end{remark}

Since the right-hand side of Eq. (\ref{main-equation}) is a discontinuous function depending on $u$, the location of interfaces between regions where $h[u]$ takes the values $+1$ and $-1$ is a priori unknown. They can be treated as the free boundaries.

\subsection{Historical review}
A first attempt to create a mathematical theory of hysteresis was made in monograph \cite{KP89} where problems with ODEs were studied. We mention also the fundamental books \cite{Vis94, BS96, Kre96} in which various hysteretic effects in spatial-distributed systems are described.

The problem (\ref{main-equation})-(\ref{initial-hysteresis}) was introduced in \cite{HJ80} where the growth of a colony of bacteria (Salmonella typhimurium) on a petri dish was modelled. The papers \cite{HJ80,HJP84} were devoted to numerical analysis of the problem, however without rigorous justification.

First existence results for solutions of problems like (\ref{main-equation})-(\ref{initial-hysteresis}) were proved in \cite{Alt85,Vis86} for modified multi-valued versions of  hysteresis operator. The main reason for these  operator extensions was the observation that the discontinuous hysteresis operator $h[u]$ is not closed with respect to topologies appropriate to its coupling with PDEs. 

In paper \cite{Alt85}, problem (\ref{main-equation})-(\ref{initial-hysteresis})   was studied in the one-(space)-dimensional case under the assumption that on the level sets $\left\lbrace u=\alpha\right\rbrace $ and $\left\lbrace u=\beta\right\rbrace $ the right-hand side of Eq. (\ref{main-equation}) allows to take any value from the whole interval $[-1,1]$. The global existence in a specially defined class of weak solutions were established there. In addition, the nonuniqueness and nonstability of such weak solutions were discussed  in \cite{Alt85} in several examples.

\begin{figure}[!h]
\centering\includegraphics[width=4in]{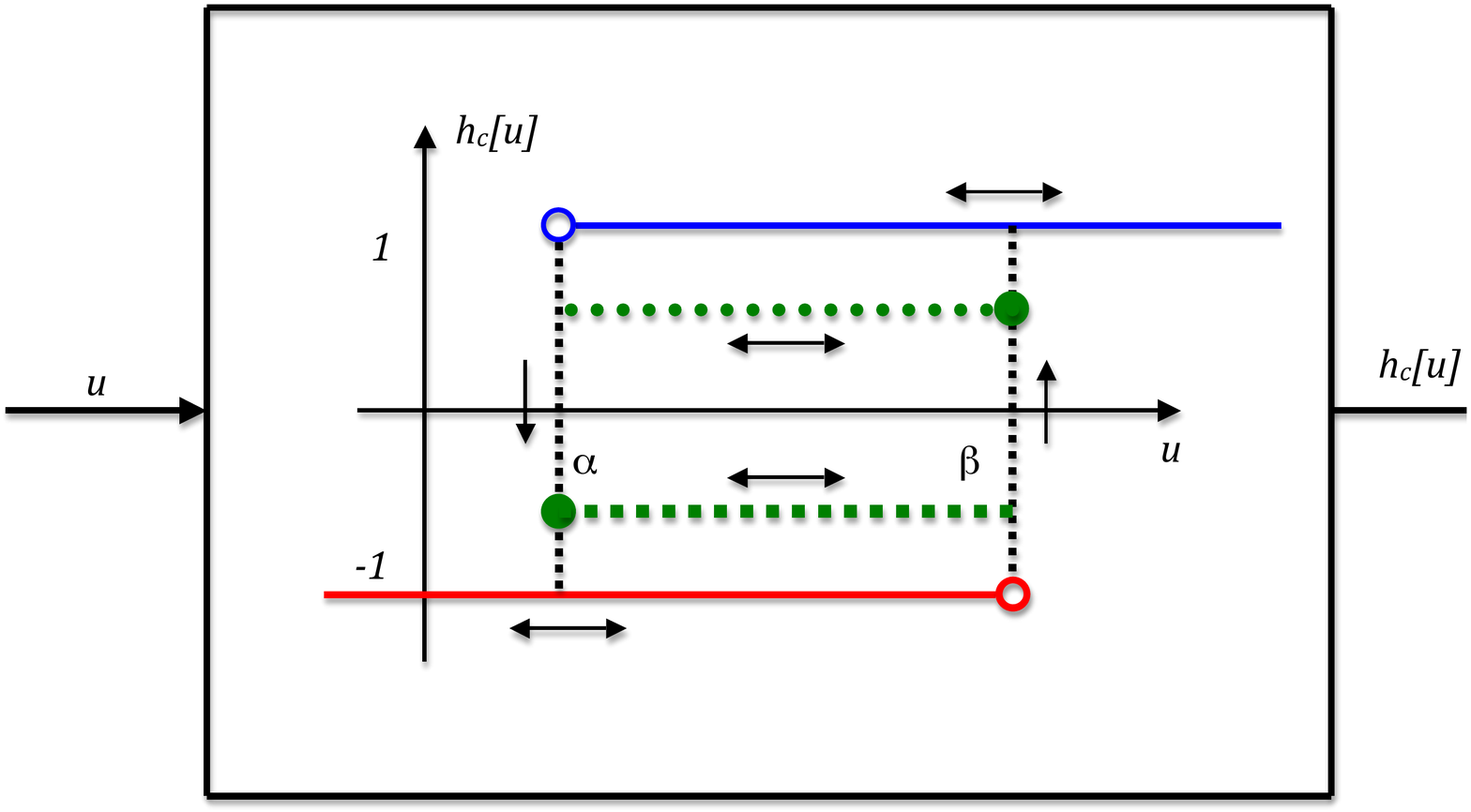}
\caption{Completed relay}
\label{fig_2}
\end{figure}

The multi-(space)-dimensional case with the non-ideal relay replaced by the so-called completed relay was treated in \cite{Vis86} where the global existence result was established for weak solutions. This completed relay operator $h_c[u]$ is the closure of $h[u]$ in suitable weak topologies. In particular, $h_c[u]$ admits any value from the whole interval $[-1,1]$ on the set $\left\lbrace \alpha \leqslant u \leqslant \beta\right\rbrace $ (see Figure~\ref{fig_2}). For  further details concerning the properties of $h_c[u]$ we refer the reader to \cite{Vis14} as well as to Chapter IV of \cite{Vis94}.

It was also shown in \cite{Vis86}  that for $\alpha=\beta$ the problem (\ref{main-equation})-(\ref{initial-hysteresis}) with $h_c[u]$ on the right-hand side of Eq. (\ref{main-equation}) is, in fact, the two-phase parabolic free boundary problem. The  properties of solutions of this two-phase problem as well as the behaviour of the corresponding free boundary  were completely studied in \cite{SUW09}.

\begin{remark}
Due to directional restrictions on jumps of the hysteresis on the tresholds $\alpha$ and $\beta$, problem  (\ref{main-equation})-(\ref{initial-hysteresis}) with $\alpha <\beta$ can not be reduced to the two-phase parabolic problem even for modified versions of the hysteresis operator.
\end{remark}

Another way to overcome the  troubles generated by discontinuity of the hysteresis operator has been proposed 
in papers 
\cite{GST13,GT12} where a \textit{special class} of strong  solutions of (\ref{main-equation})-(\ref{initial-hysteresis})  satisfying the additional \textit{transversality property}  was introduced in the one-(space)-dimensional case. This transversality property roughly speaking means that the solution $u$ has a nonvanishing spatial gradient on the free boundary. More precisely, in the case $n=1$ the transversal solutions are defined as follows.

\begin{definition} \label{def-transversal}
A function $u$ is called a  \textit{transversal solution} of (\ref{main-equation})-(\ref{initial-hysteresis}) if $u$ is a solution in the sense of Definition~\ref{def-strong-solutions} and the following hold:
\begin{enumerate}
\item if $u(\hat{x}, \hat{t}\,)=\alpha$ and $u_x(\hat{x}, \hat{t}\,)=0$ for some $\hat{x} \in \mathcal{U} \subset \mathbb{R}$ and $\hat{t}\in [0,T]$, then $h[u](x,\hat{t}\,)=-1$ in a neighborhood of $\hat{x}$;
\item if $u(\hat{x}, \hat{t}\,)=\beta$ and $u_x(\hat{x}, \hat{t}\,)=0$ for some $\hat{x} \in \mathcal{U}$ and $\hat{t}\in [0,T]$, then $h[u](x,\hat{t}\,)=1$ in a neighborhood of $\hat{x}$.
\end{enumerate}
\end{definition}
\noindent
Taking the transversal initial data and assuming that there is only a finite number of points where $\varphi$ takes the values $\alpha$ and $\beta$,   the authors of \cite{GST13} proved the local existence of strong transversal solutions of (\ref{main-equation})-(\ref{initial-hysteresis}) and showed that such solutions depend continuously on initial data. A theorem on the uniqueness of strong transversal solutions was established in \cite{GT12}.

\begin{remark}
The local well-posedness of the problem (\ref{main-equation})-(\ref{initial-hysteresis}) in the multi-(space)-dimensional case is still an open question.
\end{remark}

Recently, the approaches developed in the theory of free boundary problems have been applied to the study of strong solutions. Such an activity was initiated in 
paper \cite{AU15} where the local regularity properties of the strong solutions of  (\ref{main-equation})-(\ref{initial-hysteresis}), i.e., solutions from the Sobolev space $W^{2,1}_q$ with suffiently large $q$, were studied in the  multi-(space)-dimensional case. Without assuming the transversality property it was proved in \cite{AU15} that outside some "pathalogical" part of the free boundary we have  the optimal regularity $q=\infty$.

\begin{remark}
If the reader is not familiar with free boundary problems, we recommend him to consult the book \cite{PSU12} where the main concepts and developed methods are discussed for several model problems.
\end{remark} 

In the rest of the paper we describe in more detail the properties of the strong solutions of  (\ref{main-equation})-(\ref{initial-hysteresis}) as well as of the free boundaries that are known for today. Sections~2 and 3 are devoted to the general multi-(space)-dimensional case where we do not assume the transversality property.  Conversely, Section~4  deals with the transversal solutions in one-(space)-dimensional case. The results listed in Sections~2-3 were established in \cite{AU15}, while in Section~4 we announce the results from \cite{AU15a}.


\section{Structure of the free boundary (case $n \geqslant 1$)}
We denote 
\begin{align*}
\Omega_{\pm}(u)&=\left\lbrace z=(x,t)\in Q : \ h[u](z)=\pm 1\right\rbrace ,\\
\Gamma (u)&=\partial\Omega_+ \cap \partial\Omega_- \ \text{is the free boundary}.
\end{align*}
The latter means that $\Gamma(u)$ is the set where the function $h[u](z)$ has a jump. Note that $\Omega_{\pm}$ may consist of several components $\Omega_{\pm}^{(k)}$, respectively. In other words, 
$$
\Omega_+=\underset{k}\bigcup\,\Omega_+^{(k)}, \qquad \Omega_-=\underset{m}\bigcup\,\Omega_-^{(m)}.
$$

We also introduce special notation for some parts of $\Gamma (u)$
\begin{align*} 
\Gamma_{\alpha}(u)&:=\Gamma (u) \cap \left\lbrace u=\alpha\right\rbrace ,\\
\Gamma_{\beta} (u)&:=\Gamma (u) \cap \left\lbrace u=\beta\right\rbrace .
\end{align*}
By definition,
$$
\left\lbrace u \leqslant \alpha \right\rbrace \subset \Omega_{-}\quad \text{and}\quad \left\lbrace u \geqslant \beta\right\rbrace  \subset \Omega_{+}.
$$

Due to Definition~\ref{def-strong-solutions} it is easy to see
that the $(n+1)$-dimensional Lebesgue measure of the sets $\left\lbrace u=\alpha\right\rbrace $ and $\left\lbrace u=\beta\right\rbrace $ equals zero.

In addition, the sets $\left\lbrace u=\alpha\right\rbrace $ and $\left\lbrace u=\beta\right\rbrace $ are separated from each other. 
Moreover, for any $\varepsilon >0$ the distance from the level set $\left\lbrace u=\alpha \right\rbrace $ to the level set $\left\lbrace u=\beta\right\rbrace $ in the cylinder $Q^{\varepsilon}:=\mathcal{U}^{\varepsilon}\times (\varepsilon^2, T)$ with $\mathcal{U}^{\varepsilon} \subset \mathcal{U}$ and $\textit{dist}\left\lbrace \mathcal{U}^{\varepsilon}, \partial\mathcal{U}\right\rbrace \geqslant \varepsilon $ is estimated from below by a positive constant $\delta$  depending on $\sup\limits_{Q}|u|$, $\varepsilon$ and $\beta-\alpha$ only. 

\begin{figure}[!h]
\centering
\begin{tabular}{l|l}
\includegraphics[width=2.6in]{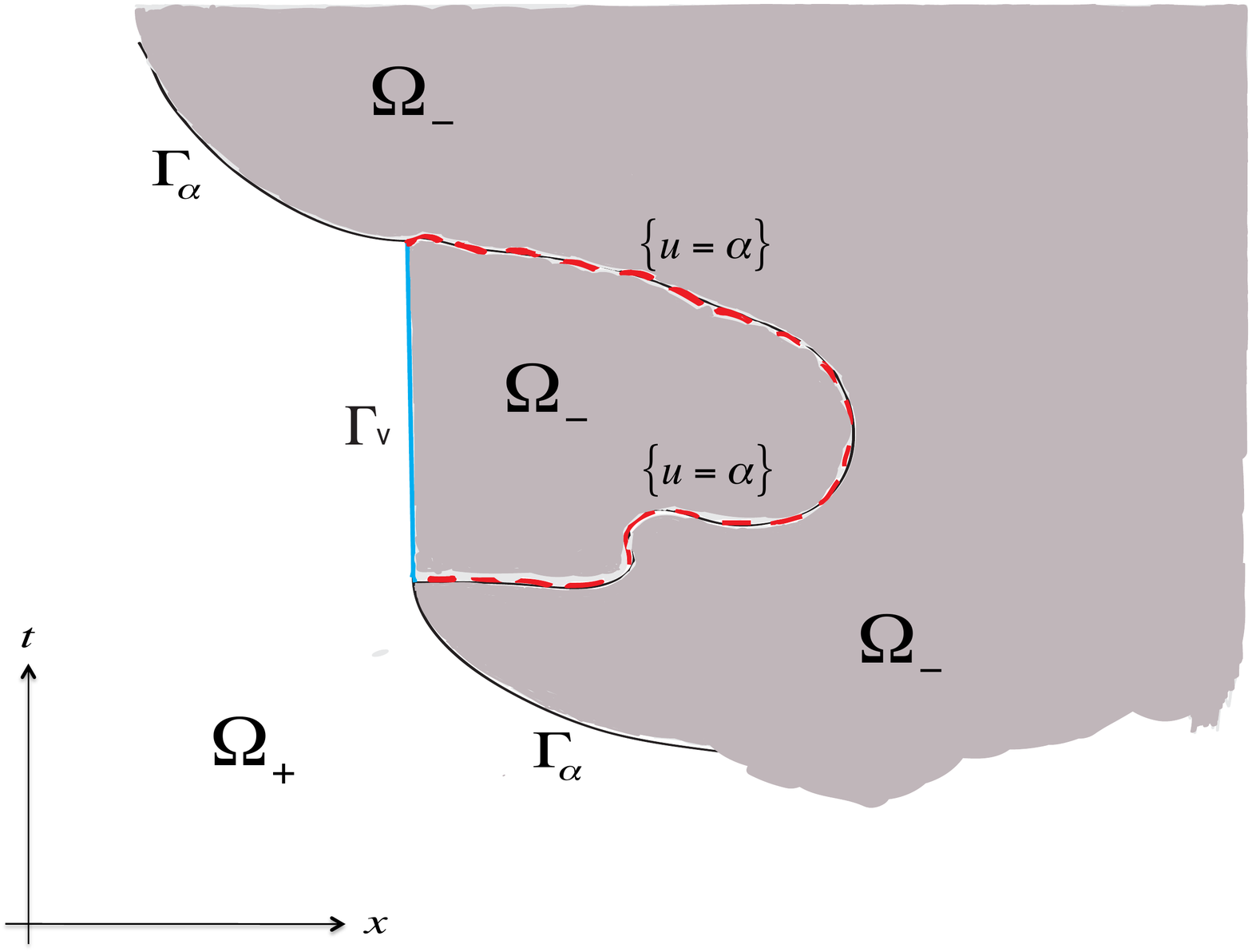} &
\includegraphics[width=2.6in]{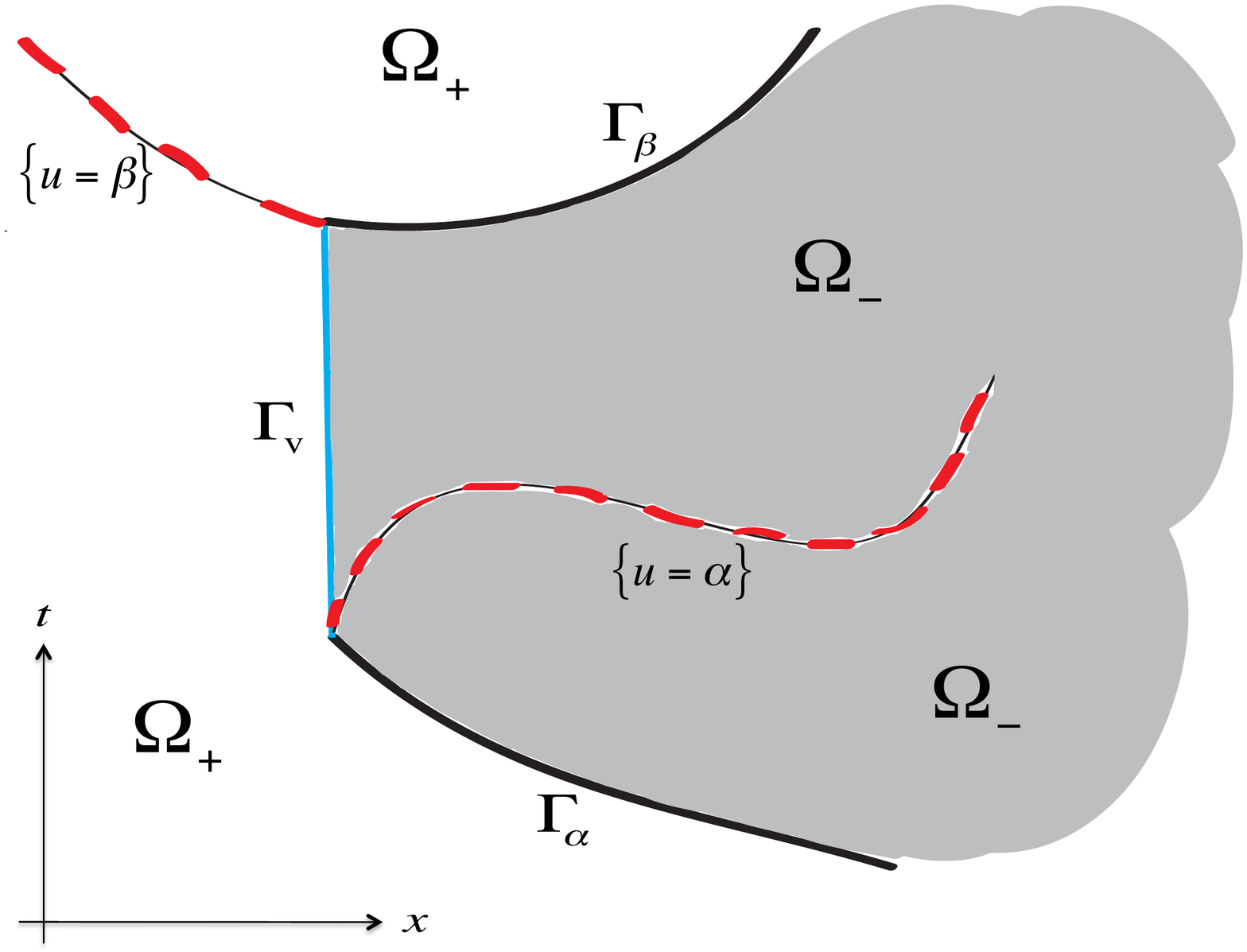}
\end{tabular}
\caption{Eventual structure of the free boundary}
\label{fig_3}
\end{figure}

Observe that the level sets $\left\lbrace u=\alpha \right\rbrace $ and $\left\lbrace u=\beta\right\rbrace $ are not alsways the parts of the free boundary $\Gamma (u)$. Indeed, if the level set $\left\lbrace u=\alpha \right\rbrace $ is locally not a $t$-graph, then a part of $\left\lbrace u=\alpha\right\rbrace $ may occur inside $\Omega_{-}$. In this case $\Gamma (u)$ 
may contain several components of $\Gamma_{\alpha}$ connected by cylindrical surfaces with generatrixes parallel to $t$-axis. Similar statement is true for the level set $\left\lbrace u=\beta\right\rbrace$.  We will denote by $\Gamma_{v}$ the set of all points $z$ lying in such  vertical parts of $\Gamma (u)$. Some examples of eventual connections of $\Gamma_v$ with $\Gamma_{\alpha}$ and $\Gamma_{\beta}$ are shown on Figure~\ref{fig_3}. 

Recall that by definition of $\Gamma_{\alpha}(u)$ the function $h[u]$ has a jump in $t$-direction from $+1$ to $-1$ there. The latter means that if we cross the free boundary $\Gamma_{\alpha}(u) \setminus \Gamma_v$ in positive $t$-direction then the corresponding phases change from $\Omega_{+}$ to $\Omega_{-}$ (see Figure~\ref{fig_3} again). The similar statement can be made for the neighborhood of $\Gamma_{\beta} (u)$, that is if we cross the free boundary $\Gamma_{\beta}(u)\setminus \Gamma_v$ in positive $t$-direction then the phases will change from $\Omega_-$ to $\Omega_+$.

Thus, we have
$$
\Gamma (u)=\Gamma_{\alpha} (u) \cup \Gamma_{\beta} (u) \cup \Gamma_{v}.
$$

It should be noted that this $\Gamma_v$ is just the "pathalogical" part of the free boundary  mentioned in Introduction. Indeed, we have no information about the values of $u$ on $\Gamma_v$, since $\Gamma_v$ is, in general, not the level set $\left\lbrace u=\alpha\right\rbrace $ as well as not the level set $\left\lbrace u=\beta\right\rbrace $. Furthermore, for any direction $e\in \mathbb{R}^n$ functions $\left( D_eu\right)_+= \max\left\lbrace D_eu,0\right\rbrace $ and $\left( D_eu\right)_-= \max\left\lbrace -D_eu,0\right\rbrace $ are, in general, not sub-caloric near $\Gamma_v$. The latter fact causes the serious difficulties in studying the regularity properties of solutions.

We will also distinguish the following parts of $\Gamma$:
$$
\Gamma_{\alpha}^0 (u)=\Gamma_{\alpha} (u) \cap \left\lbrace |Du|=0\right\rbrace , \qquad \Gamma_{\alpha}^* (u)=\Gamma_{\alpha} (u) \setminus \Gamma_{\alpha}^0 (u).
$$
The sets $\Gamma_{\beta}^0$ and $\Gamma^*_{\beta}$ are defined analogously. In addition, we set
$$
\Gamma^0(u):=\Gamma^0_{\alpha}(u) \cup \Gamma^0_{\beta} (u), \quad
\Gamma^*(u):=\Gamma^*_{\alpha}(u) \cup \Gamma^*_{\beta} (u).
$$

Using the von Mises transformation combined with the parabolic theory and the Implicit Function Theorem it is possible to show that $\Gamma^* (u)\setminus \Gamma_{v}$ is locally a $C^1$-surface.

\section{Optimal~regularity~of~$u$~beyond~$\Gamma_{v}$~(case $n~\geqslant~1$)}
In this section we discuss the\textit{ local} regularity properties of a strong solution $u$ of (\ref{main-equation})-(\ref{initial-hysteresis}). 

For simplicity of notation we will denote by $M$ a majorant of $\sup_{Q}|u|$ and will highlight below the dependence of all obtained estimates on $M$.

Recall that 
the general parabolic theory (see, e.g. \cite{LSU67}) provides for any $\varepsilon>0$ the estimates
$$
\|\partial_t u\|_{q, Q^{\varepsilon}}+\|D^2u\|_{q, Q^{\varepsilon}}\leqslant N_1(\varepsilon, q,M) \quad \forall q<\infty, 
$$
where $Q^{\varepsilon}=\mathcal{U}^{\varepsilon} \times (\varepsilon^2, T)$, $\mathcal{U}^{\varepsilon} \subset \mathcal{U}$ and $\textit{dist}\, \left\lbrace \mathcal{U}^{\varepsilon}, \partial\mathcal{U}\right\rbrace \geqslant \varepsilon$. 
Moreover, the function $u$ and its spatial gradient $Du$ are H{\"o}lder continuous in $Q$, and $u\in C^{\infty}$ in the interior of the sets $\Omega_{\pm}$, as well.

We note also that 
if  $\partial\mathcal{U}$ as well as the values of $u$ on the parabolic boundary of $Q$ are smooth
then  the corresponding estimates of $L^q$-norm for $\partial_t u$ and $D^2u$ are true  in the whole cylinder $Q$.

Contrariwise, $\Delta u-\partial_t u$ has a jump across the free boundary $\Gamma (u)$. Thus, $W^{2,1}_{\infty, loc}$ is the best possible regularity of solutions.

The optimal (i.e., $W^{2,1}_{\infty}$) regularity is not obvious. A crucial point here is  the quadratic growth estimate of the type
\begin{equation} \label{quadratic-growth-u}
\sup\limits_{Q_r(z^0)}|u| \leqslant N_2(\rho_0,\epsilon, M)r^2 \quad \text{for} \quad r\leqslant \min \left\lbrace \rho_0, \epsilon \right\rbrace 
\end{equation}
with $z^0=(x^0,t^0)\in \Gamma^0(u)\setminus \Gamma_v$. Here $\rho_0$ denotes  the parabolic distance from $z^0$ to $\Gamma_v$, while $\epsilon$ stands for the parabolic distance from $z^0$ to $\partial' Q$ and $$Q_r(z^0):=\left\lbrace x : |x-x^0|<r\right\rbrace \times (t^0-r^2, t^0+r^2).
$$
\begin{remark}
The parabolic distance $\textit{dist}_p$ from a point $z^0=(x^0,t^0)$ to a set $\mathcal{D} \subset \mathbb{R}^{n+1}$ is defined as 
$$
\textit{dist}_p\left\lbrace z^0, \mathcal{D}\right\rbrace :=\sup\limits_{\rho>0}\left\lbrace  Q_{\rho}(z^0)\cap \left\lbrace t \leqslant t^0\right\rbrace  \cap
\mathcal{D}=\emptyset \right\rbrace. 
$$
\end{remark}
To show  quadratic bound (\ref{quadratic-growth-u}) we argue by a contradiction and combine this with a local rescaled version of the famous Caffarelli monotonicity formula.
\begin{remark}
More detailed information about Caffarelli's monotonicity formula and its local rescaled version can be found in  book \cite{CS05} and in papers \cite{AU15, ASU00}, respectively.
\end{remark}
Further, it can be verified
that quadratic growth estimate (\ref{quadratic-growth-u}) implies the corresponding linear bound for $|Du|$
\begin{equation} \label{linear-growth-Du}
\sup\limits_{Q_r(z^0)}|Du| \leqslant N_2(\rho_0,\epsilon, M)r \quad \text{for all}\quad r \leqslant \min \left\lbrace \rho_0, \epsilon\right\rbrace 
\end{equation}
with $z^0 \in \Gamma^0(u) \setminus \Gamma_v$.

The dependence of $N_2$ on the distance $\rho_0$ in (\ref{quadratic-growth-u})-(\ref{linear-growth-Du}) arises due to the monotonicity formula. Unfortunately, near $\Gamma_v$ neither  the local rescaled version of Caffarelli's monotonicity formula nor its generalisation such as the almost monotonicity formula introduced in \cite{EP08} are applicable to positive and negative parts of the space directional derivatives $D_eu$. 

Besides estimates (\ref{quadratic-growth-u})-(\ref{linear-growth-Du}),    the information about
behaviour of $\partial_tu$ near $\Gamma^* (u)$ plays an important role. As already mentioned above, $\partial_t u$ may has jumps across the free boundary. Actually,
one can show that $\partial_t u$ is a continuous function in a neighborhood of $z^* \in \Gamma^* (u) \setminus~\Gamma_v$. In addition, the monotonicity of jumps of $h[u]$ in $t$-direction provides  the one-sided estimates of the time derivative of $u$ near $\Gamma (u)$. More precisely, the estimate of $\partial_t u$ from below holds true  near $\Gamma_{\alpha}(u)$, whereas the estimate of $\partial_tu$ from above holds true near $\Gamma_{\beta} (u)$. Combination of these  results  with the observation that $\partial_tu \leqslant 0$ on $\Gamma_{\alpha}^* (u)\setminus \Gamma_v$ and $\partial_t u \geqslant 0$ on $\Gamma_{\beta}^* (u)\setminus \Gamma_v$ leads to an absolute estimate of the time-derivative of $u$ on the set $\Gamma^*~(u)~\setminus~\Gamma_v$. Namely,
\begin{equation} \label{time-derivative-estimate}
\sup_{\Gamma^*(u)\setminus \Gamma_v}|\partial_tu| \leqslant N_3(\epsilon, M, \beta-\alpha),
\end{equation}
where the constant $N_3$, contrary to $N_2$, depends only on given quantities. In addition, we can show that the mixed second derivatives $D_i(\partial_t u)$ are $L^2_{loc}$-functions in $Q \setminus \left( \Gamma^0(u)\cup \Gamma_v\right) $.

Now, estimates (\ref{quadratic-growth-u})-(\ref{time-derivative-estimate}) allow us to apply the methods from the theory of free boundary problems and  estimate $|\partial_tu(z)|$ and $|D^2u(z)|$ for any $z$ being a point of smothness for $u$.
The corresponding  details can be found in \cite{AU15}. The main result is formulated as follows.
\begin{theorem} 
Let $u$ be a solution of Eq.~\ref{main-equation}, and let $z\in Q\setminus \Gamma (u)$. Then
\begin{equation} \label{optimal-regularity}
|\partial_t u(z)|+|D^2u(z)| \leqslant C(\rho_0, \epsilon, M, \beta-\alpha).
\end{equation}
Here $\rho_0:=\textit{dist}_p\left\lbrace z, \Gamma_v\right\rbrace $ and $\epsilon:=\textit{dist}_p\left\lbrace z, \partial'Q\right\rbrace $. 
\end{theorem}

We emphasize that the constant $C$ in (\ref{optimal-regularity}) does not depend on the parabolic distance from $z$ to $\Gamma^0(u)$ as well as to $\Gamma^*(u)$. Unfortunately, we cannot remove the dependence of $C$ on $\rho_0$, i.e. on the parabolic distance of $z$ to $\Gamma_v$.


\section{Transversal solutions (case $n=1$)}
For one-(space)-dimensional  transversal solutions of (\ref{main-equation})-(\ref{initial-hysteresis}) the results of Sections 2 and 3 can be strengthened.

As in the papers \cite{GST13, GT12} we will suppose that initially there are only finite number of different components $\Omega_{\pm}^{(k)}$.


It is obvious that for $n=1$ the "pathalogical" part of the free boundary $\Gamma_v$ is a union of vertical segments parallel to $t$-axis.
Further, according to Definition~\ref{def-transversal} we conclude that for transversal solutions $u$ the inequality
$$
|u_x(x,t)|>0 
$$
holds true on $\Gamma_{\alpha}(u) \cup \Gamma_{\beta} (u)$. However, on $\Gamma_v$ the function $u_x$ may vanish even for transversal solutions. In addition, the time-derivative $\partial_t u$ is continuous across $\Gamma_v$ except  eventually the end points of the vertical segments.
Thus,
for transversal solutions we have $\Gamma^0(u)=\emptyset$, and, consequently,  $\Gamma^*_{\alpha}(u)=\Gamma_{\alpha}(u)$ and $\Gamma^*_{\beta}(u)~=~\Gamma_{\beta}(u)$.

\begin{figure}[!h]
\centering\includegraphics[width=4in]{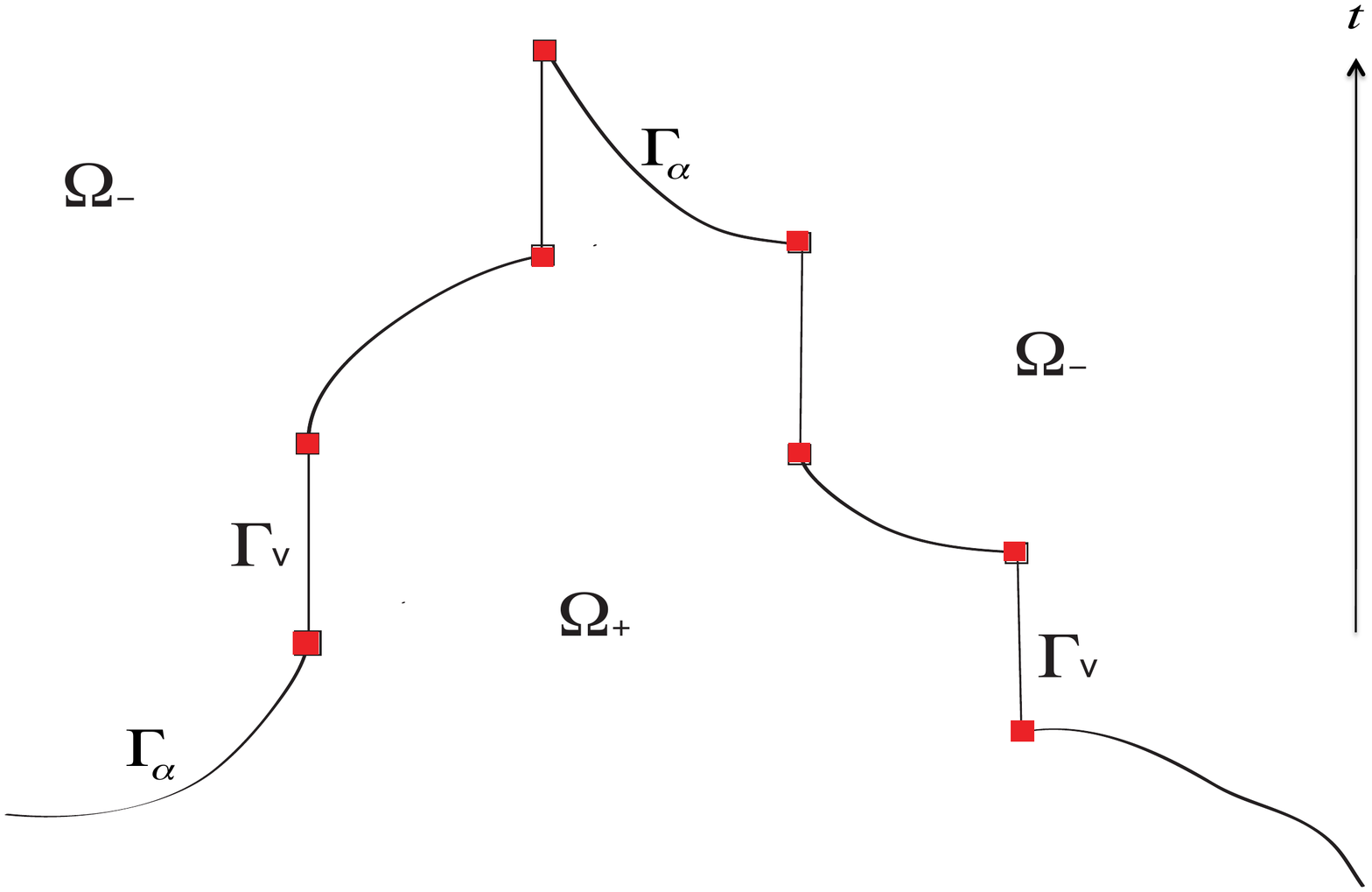}
\caption{Local behaviour of $\Gamma (u)$}
\label{fig_4}
\end{figure}

The transversality property of a solution $u$  implies the monotone behaviour of the free boundary $\Gamma (u)$. Indeed, it is possible to show that   the set
$\Omega_+^{(k)}$ is locally either a  subgraph of a monotone curve or a subgraph of a union of two monotone curves  with different character  of monotonicity. These monotone curves include the parts of $\Gamma_{\alpha}(u)$ and also may contain vertical segments. An example of a possible union of  two monotone curves  is provided on Figure~\ref{fig_4}.








Similar statement is true for the set $\Omega^{(k)}_{-}$ near $\Gamma_{\beta}(u)$. So, $\Omega^{(k)}_{-}$
is locally either a  subgraph of a monotone curve  or a subgraph of a union of two monotone curves with different character  of monotonicity.

Finally, under the additional assumption that on a time-interval $(t^1,t^2)$ the free boundary part $\Gamma_v$ consists of at most finite number of vertical segments, it is possible to prove the optimal regularity result  up to $\Gamma_v$.
In other words,
one can establish that  for any $z\in Q\cap (t^1_i,t^2_i)$  being a point of smoothness for transversal solutions $u$ the constant $C$ from (\ref{optimal-regularity}) is independent on the parabolic distance from $z$ to $\Gamma_v$. 



The detailed explanation of all results presented in this section can be found in~\cite{AU15a}.



\section{Conclusion}
Below we discuss several open problems and areas of further work. \vspace{0.2cm}

It should be emphasized that quadratic growth estimate (\ref{quadratic-growth-u}) allows to apply the standard parabolic scaling in points $z^0\in \Gamma^0 (u)\setminus \Gamma_v$  and to obtain the corresponding  blow-up limits~of~$u$. Using a version of parabolic Weiss's monotonicity formula  it can be shown that all blow-ups are homogeneous functions. Classification of all possible blow-up limits and further studying the regularity of the free boundary $\Gamma (u) \setminus \Gamma_v$ are open problems. For $n=1$ some preliminary results of this kind are established in \cite{AU15a}.

\begin{remark}
The original version of parabolic Weiss's monotonicity formula can be found in \cite{W99}.
\end{remark}

Let $u$ satisfy on $\partial'Q$ the inequalities $\alpha \leqslant u \leqslant \beta$. In this case even for weak solutions we have the same inequalities $\alpha \leqslant u \leqslant \beta$ a.e. in $Q$. We suppose that strong solutions of (\ref{main-equation})-(\ref{initial-hysteresis}) do not exist in this case. \vspace{0.2cm}

Note that numerical examples given in presentation \cite{GTR} suggest  the possible non-existence of the strong solutions with nontransveral initial conditions. Confirmation or  rejection of this hypothesis is an open question even for the one-(space)-dimensional case. \vspace{0.2cm}

Another challenging problem is to define  the transversality property for strong multi-(space)-dimensional solutions. Some results in this direction were obtained very recently in \cite{C14}. \vspace{0.2cm}

The last (but not least) hypothesis concerns the regularity of the free boundary in the nontransversal case $n=1$. It is proposed to prove or disprove the assertion that the free boundary $\Gamma(u)$ is smooth except the end-points of the vertical segments provided $\Gamma_v$ consists of at most finite number of such segments.
\vspace{0.2cm}

\section*{Acknowledgment}

This work was supported by the Russian Foundation of Basic Research (RFBR) through the grant number
14-01-00534,  by the St. Petersburg State University grant 6.38.670.2013 and by the grant "Nauchnye Shkoly", NSh-1771.2014.

The authors also thank the Isaac Newton Institute for Mathematical Sciences, Cambridge, UK, where a part of this work was done during the program \textit{Free Boundary Problems and Related Topics}.



\Addresses

\begin{thebibliography}{9}

\bibitem{Vis86} Visintin A. 1986 Evolution problems with hysteresis in the source term. \textit{SIAM J. Math. Anal.} \textbf{17}(5):
1113-1138.

\bibitem{Alt85} Alt H.W. 1985 On the thermostat problem. \textit{Control Cybernet.} \textbf{14}(1-3): 171-193.

\bibitem{HJ80} Hoppensteadt F.C., J{\"a}ger W. 1980 Pattern formation by bacteria. In \textit{Biological growth and spread (Proc. Conf., Heidelberg, 1979)}, volume 38 of \textit{Lecture Notes in Biomath.}, pp. 68-81. Berlin-New-York: Springer.

\bibitem{HJP84}  Hoppensteadt F.C., J{\"a}ger W., P{\"o}ppe C. 1984  A hysteresis model for bacterial growth patterns. In \textit{Modellimg of patterns in space and time (Heidelberg, 1983)}, volume 55 of \textit{Lecture Notes in Biomath.}, pp. 123-134. Berlin: Springer. 

\bibitem{K06} Kopfova J. 2006 Hysteresis and biological models.  \textit{J. Phys. Conference Series.} \textbf{55}: 130-134.

\bibitem{KP89} Krasnosel'ski\u{i} M.A., Pokrovski\u{i} A.V. 1989 \textit{Systems with Hysteresis.} (Translated from Russian: "Sistemy s gisterezisom", Nauka, Moscow, 1983). Springer-Verlag. Berlin.

\bibitem{Vis94} Visintin A. 1994 \textit{Differential models of hysteresis. Applied Mathematical Sciences, 111.} Springer-Verlag. Berlin.

\bibitem{BS96} Brokate M., Sprekels J. 1996 \textit{Hysteresis and phase transitions. Applied Mathematical Sciences, 121.} Springer-Verlag. New York.

\bibitem{Kre96} Krej\v{c}\'{i} P. 1996 \textit{Hysteresis, convexity and dissipation.} GAKUTO International Series. Mathematical Sciences and Applications, 8. Gakk\={o}tosho Co., Ltd., Tokyo.

\bibitem{Vis14} Visintin A. 2014 Ten issues about hysteresis. \textit{Acta Appl. Math.} \textbf{132}: 635-647.

\bibitem{SUW09} Shahgholian H., Uraltseva N., Weiss G.S. 2009 A parabolic two-phase obstacle-like equation. \textit{Adv. Math.} \textbf{221}(3): 861-881.

\bibitem{GST13} Gurevich P., Shamin R., Tikhomirov S. 2013 Reaction-diffusion equations with spatially distributed hysteresis.
\textit{SIAM J. Math. Anal.} \textbf{45}(3): 1328-1355.

\bibitem{GT12} Gurevich P., Tikhomirov S. 2012 Uniqueness of transverse solutions for reaction-diffusion equations with spatially distributed hysteresis. \textit{Nonlinear Anal.} \textbf{75}(18): 6610-6619.

\bibitem{AU15} Apushkinskaya D.E., Uraltseva N.N. 2015 On regularity properties of solutions to the hysteresis-type problems. \textit{Interfaces and Free~Boundaries}~\textbf{17}.

\bibitem{AU15a} Apushkinskaya D.E., Uraltseva N.N. Uniform estimates of transversal solutions to the one-space-dimensional hysteresis-type problem. \textit{In preparation}.

\bibitem{PSU12} Petrosyan A., Shahgholian H., Uraltseva N. 2012 \textit{Regularity of free boundaries in obstacle type problems. Graduate Studies in Mathematics, 136.} American Mathematical Society. Providence, RI.

\bibitem{LSU67} Lady\v{z}enskaja O.A., Solonnikov V.A., Ural'ceva N.N. 1967 \textit{Linear and quasilinear equations of parabolic type.} Translations of Mathematical Monographs, Vol.23. American Mathematical Society. Providence, RI.

\bibitem{CS05} Caffarelli L., Salsa S. 2005 \textit{A geometric approach to free boundary problems. Graduate Studies in Mathematics, 68.} American Mathematical Society. Providence, RI.

\bibitem{ASU00} Apushkinskaya D.E., Shahgholian H., Uraltseva N.N. 2000 Boundary estimates for solutions of a parabolic free boundary problem. \textit{Zap. Nauchn. Sem. S.-Petersburg. Otdel. Mat. Inst. Steklov. (POMI)} \textbf{271}: 39-55.

\bibitem{EP08} Edquist A., Petrosyan A. 2008 A parabolic almost monotonicity formula. \textit{Math. Ann.} \textbf{341}(2): 429-454.

\bibitem{W99} Weiss G.S. 1999 Self-similar blow-up and Hausdorff dimension estimates for a class of parabolic free boundary problems. \textit{SIAM J. Mat. Anal.} \textbf{30}(3):623-644.



\bibitem{GTR} Gurevich P., Tikhomirov S., Ron E.  Reaction-diffusion equations with discontinuous hysteresis (poster). 
See https://sites.google.com/site/sergeytikhomirov/

\bibitem{C14} Curran M. 2014 Local well-poseness of a reaction-diffusion equation with hysteresis. Master thesis, Fachbereich Mathematik und Informatik, Freie Universit{\"a}t Berlin.





\end{thebibliography}
\end{document}